\documentclass[12pt]{amsart}
\usepackage{amssymb}
\usepackage{amsfonts}
\usepackage{amsbsy}
\usepackage{latexsym}
\newtheorem{theorem}{Theorem}[section]
\newtheorem{lemma}[theorem]{Lemma}
\newtheorem{corollary}[theorem]{Corollary}

\theoremstyle{definition}
\newtheorem{definition}[theorem]{Definition}

\begin{document}
\title {Additivity of Maps on Triangular algebras}
\author{Xuehan Cheng}
\address{College of  Mathematics and Information, Ludong University,  Yantai, 264025 P. R. China}

\author{Wu Jing}
\address{Department of Mathematics and Computer Science, Fayetteville State University,  Fayetteville, NC 28301}
 \email{wjing@uncfsu.edu}

\thanks{This work is partially supported by NNSF of China (No. 10671086) and NSF of Ludong University (No. LY20062704).}
\date{June 8, 2007}
\begin{abstract}
We prove that every multiplicative bijective map, Jordan bijective
map and Jordan triple bijective map from a triangular algebra onto
any ring is automatically additive. \end{abstract} \maketitle

\section{Introduction}

If a ring $\mathcal {R}$ contains a nontrivial idempotent, it is
kind of surprise that  every multiplicative bijective map from
$\mathcal {R}$ onto an arbitrary ring is automatically additive.
This result was given by Martindale III in his excellent paper \cite
{ma695}. More precisely, he proved:
\begin{theorem}(\cite{ma695}) Let $\mathcal {R}$ be a ring containing a family $\{ e_{\alpha }:\alpha \in \Lambda \} $ of idempotents which satisfies:

(i) $x\mathcal {R}=\{ 0\} $ implies $x=0$;

(ii) If $e_{\alpha }\mathcal {R}x=\{ 0\} $ for each $\alpha \in
\Lambda $, then $x=0$ (and hence $\mathcal {R}x=\{ 0\} $ implies
$x=0$);

(iii) For each $\alpha \in \Lambda $, $e_{\alpha }xe_{\alpha
}\mathcal {R}(1-e_{\alpha })=\{ 0\} $ implies $e_{\alpha }xe_{\alpha
}=0$.

Then any multiplicative bijective map from $\mathcal {R}$ onto an
arbitrary ring $\mathcal {R}^{\prime }$ is additive.
\end{theorem}

Note that the proof of \cite {ma695} has become a  standard argument
and been applied widely by several authors in investigating the
additivity of  maps on rings as well as  on operator algebras (see
\cite {ji}-\cite{lu2273}). Following this standard argument (see
\cite {ma695}), in this paper we continue to study the additivity of
map on triangular algebras.  We will show that every multiplicative
bijective  map, Jordan bijective map, and Jordan triple bijective
map from a triangular algebra onto an arbitrary ring is additive.
The advantage of triangular algebra is that its special structure
enables us to borrow the idea from \cite {ma695} while we do not
require the existence of  nontrivial idempotents.

We now  introduce some definitions and results.

\begin{definition} Let $\mathcal {R}$ and $\mathcal {R}^{\prime }$ be two rings, and $\phi $ be a map from $\mathcal {R}$ to $\mathcal {R}^{\prime }$. Suppose that  $a, b$, and $c$ are arbitrary elements of $\mathcal {R}$.

$\phi $ is said to be \textit{multiplicative} if
$$\phi (ab)=\phi (a)\phi (b).$$

$\phi $ is called a \textit{Jordan map} if
$$\phi (ab+ba)=\phi (a)\phi (b)+\phi (b)\phi (a).$$

Finally, $\phi $ is called a \textit{Jordan triple map} if
$$\phi (abc+cba)=\phi (a)\phi (b)\phi (c)+\phi (c)\phi (b)\phi (a).$$
\end{definition}
Note that no additivity  is required in this definition.

 Recall that a \textit{triangular algebra} $Tri(\mathcal {A}, \mathcal {M}, \mathcal {B})$ is  an algebra of the form

 \begin{displaymath}Tri(\mathcal {A}, \mathcal {M}, \mathcal {B})=\{
\left( \begin{array}{ll}
a & m\\
0 & b
\end{array}\right): a\in \mathcal {A}, m\in \mathcal {M}, b\in \mathcal {B}\}
\end{displaymath}
under the usual matrix operations, where $\mathcal {A}$ and
$\mathcal {B}$ are two algebras over a commutative ring $\mathcal
{R}$, and $\mathcal {M}$ is an $(\mathcal {A}, \mathcal
{B})$-bimodule which is faithful as a left $\mathcal {A}$-module and
also as a right $\mathcal {B}$-module (see \cite{cheung}).

Throughout this paper, we set \begin{displaymath} \mathcal
{T}_{11}=\{ \left( \begin{array}{ll}
a & 0\\
0 & 0
\end{array}\right): a\in \mathcal {A} \},
\end{displaymath}
\begin{displaymath}\mathcal {T}_{12}=\{
\left( \begin{array}{ll}
0 & m\\
0 & 0
\end{array}\right): m\in \mathcal {M} \},
\end{displaymath}
and
\begin{displaymath} \mathcal {T}_{22}=\{
\left( \begin{array}{ll}
0 & 0\\
0 & b
\end{array}\right): b\in \mathcal {B} \}.
\end{displaymath}
Then we may write $\mathcal {T}=\mathcal {T}_{11}\oplus \mathcal
{T}_{12}\oplus \mathcal {T}_{22}$, and every element $a\in\mathcal
{T}$ can be written as $a=a_{11}+a_{12}+a_{22}$. Note that notation
$a_{ij}$ denotes an arbitrary element of
 $\mathcal {T}_{ij}$.

Let $X$ be a Banach space. We denote by $B(X)$ the algebra of all
bounded linear operators on $X$. A subalgebra $\mathcal {A}$ of
$B(X)$ is called a \textit{standard operator algebra} if $\mathcal
{A}$ contains all finite rank operator. Note that if $A\in \mathcal
{A}$ and $A\mathcal {A}=\{ 0\} $ (or $\mathcal {A}A=\{ 0\} $), then
$A=0$.

\section{Multiplicative Maps on Triangular Algebras}
In this section we study the additivity of multiplicative map on
triangular algebras. We now state our first main result.

\begin{theorem} \label{theorem}   Let $\mathcal {A}$ and $\mathcal {B}$ be two algebras over a commutative ring $\mathcal {R}$, $\mathcal {M}$ a faithful $(\mathcal {A}, \mathcal {B})$-bimodule, and $\mathcal {T}$ be the triangular algebra $Tri(\mathcal {A}, \mathcal {M}, \mathcal {B})$. Suppose that algebras $\mathcal {A}$ and $\mathcal {B}$ satisfy:

(i) If $a\mathcal {A}=\{ 0\} $, or $\mathcal {A}a=\{ 0\} $, then
$a=0$;

(ii) If $b\mathcal {B}=\{ 0\} $, or $\mathcal {B}b=\{ 0\} $, then
$b=0$.

Then any multiplicative bijective map from $\mathcal {T}$ onto an
arbitrary ring $\mathcal {R}^{\prime }$ is additive.
\end{theorem}

 The proof of  this theorem  is organized into  a series of lemmas. In what follows, $\phi $ will be a multiplicative  bijective map from $\mathcal {T}$ onto an arbitrary ring $\mathcal {R}^{\prime }$.

\begin{lemma} $\phi (0)=0$.
\end{lemma}
\begin{proof} Since $\phi $ is surjective, there exists $a\in \mathcal {T}$ such that $\phi (a)=0$. Then $\phi (0)=\phi (0\cdot a)=\phi (0)\phi (a)=\phi (0)\cdot 0=0$.
\end{proof}

\begin{lemma}\label{1112}For any $a_{11}\in \mathcal {T}_{11}$ and $b_{12}\in \mathcal {T}_{12}$, we have
$$\phi (a_{11}+b_{12})=\phi (a_{11})+\phi (b_{12}).$$
\end{lemma}
\begin{proof}
Let $c\in \mathcal {T}$ be chosen such that $\phi (c)=\phi
(a_{11})+\phi (b_{12})$.

For arbitrary $t_{11}\in \mathcal {T}_{11}$, we consider
\begin{eqnarray*}
\phi (ct_{11})&=&\phi (c)\phi (t_{11})=(\phi (a_{11})+\phi (b_{12}))\phi (t_{11})\\
&=&\phi (a_{11})\phi (t_{11})+\phi (b_{12})\phi (t_{11})=\phi
(a_{11}t_{11}).
\end{eqnarray*}
Hence, $ct_{11}=a_{11}t_{11}$, and so $c_{11}=a_{11}$.

Similarly, we can get $c_{22}=0$.

We now show that $c_{12}=b_{12}$. For any $t_{11}\in \mathcal
{T}_{11}$ and $s_{22}\in \mathcal {T}_{22}$, we obtain
\begin{eqnarray*}
\phi (t_{11}cs_{22})&=&\phi (t_{11})\phi (c)\phi (s_{22})\\
&=&\phi (t_{11})(\phi (a_{11}+\phi (b_{12}))\phi (s_{22})\\
&=&\phi (t_{11}b_{12}s_{22}).
\end{eqnarray*}
It follows that $t_{11}cs_{22}=t_{11}b_{12}s_{22}$, which give us
that $c_{12}=b_{12}$.
\end{proof}

Similarly, we have the following lemma.
\begin{lemma} \label{2212}For arbitrary $a_{22}\in \mathcal {T}_{22}$ and $b_{12}\in \mathcal {T}_{12}$, the following holds true.
$$\phi (a_{22}+b_{12})=\phi (a_{22})+\phi (b_{12}).$$
\end{lemma}

\begin{lemma}\label{11121222}
For any $a_{11}\in \mathcal {T}_{11}$, $b_{12}, c_{12}\in \mathcal
{T}_{12}$, and $d_{22}\in \mathcal {T}_{22}$, we have
$$\phi (a_{11}b_{12}+c_{12}d_{22})=\phi (a_{11}b_{12})+\phi (c_{12}d_{22}).$$
\end{lemma}
\begin{proof} By Lemma \ref{1112} and Lemma \ref{2212}, we have
\begin{eqnarray*}
& &\phi (a_{11}b_{12}+c_{12}d_{22})\\
&=&\phi (a_{11}+c_{12})(b_{12}+d_{22})\\
&=&\phi (a_{11}+c_{12})\phi (b_{12}+d_{22})\\
&=&(\phi (a_{11})+\phi (c_{12}))(\phi (b_{12})+\phi (d_{22}))\\
&=&\phi (a_{11})\phi (b_{12})+\phi (a_{11})\phi (d_{22})+\phi (c_{12})\phi (b_{12})+\phi (c_{12})\phi (d_{22})\\
&=&\phi (a_{11}b_{12})+\phi (c_{12}d_{22}).
\end{eqnarray*}
\end{proof}
\begin{lemma}\label{12} $\phi $ is additive on $\mathcal {T}_{12}$.
\end{lemma}
\begin{proof} Suppose that $a_{12}$ and $b_{12}$ are two elements of $\mathcal {T}$.
We pick $c\in \mathcal {T}$ such that $\phi (c)=\phi (a_{12})+\phi (b_{12})$.
For ant $t_{11}\in \mathcal {T}_{11}$ and $s_{22}\in \mathcal
{T}_{22}$, we compute
\begin{eqnarray*}
\phi (t_{11}cs_{22})&=&\phi (t_{11})\phi (c)\phi (s_{22})\\
&=&\phi (t_{11})(\phi (a_{12}+\phi (b_{12}))\phi (s_{22})\\
&=&\phi (t_{11})\phi (a_{12})\phi (s_{22})+\phi (t_{11})\phi (b_{12})\phi (s_{22})\\
&=&\phi (t_{11}a_{12}s_{22}+t_{11}b_{12}s_{22}).
\end{eqnarray*}
Note that in the last equality we apply Lemma \ref{11121222}. It
follows that $t_{11}cs_{22}=t_{11}a_{12}s_{22}+t_{11}b_{12}s_{22}$,
and so $c_{12}=a_{12}+b_{12}$.

We can get $c_{11}=c_{22}=0$ by considering $\phi (ct_{11})$ and
$\phi (t_{22}c)$ for arbitrary $t_{11}\in \mathcal {T}_{11}$ and
$t_{22}\in \mathcal {T}_{22}$ respectively.
\end{proof}

\begin{lemma}\label{11}
$\phi $ is additive on $\mathcal {T}_{11}$.
\end{lemma}
\begin{proof}
Let $a_{11}, b_{11}\in \mathcal {T}_{11}$, and $c\in \mathcal {T}$
such that $\phi (c)=\phi (a_{11}+\phi (b_{11})$. We only show that
$c_{11}=a_{11}+b_{11}$. One can easily get $c_{22}=c_{12}=0$ by
considering $\phi (t_{22}c)$ and $\phi (t_{11}cs_{22})$ for any
$t_{11}\in \mathcal {T}_{11}$ and $t_{22}, s_{22}\in \mathcal
{T}_{22}$.

For ant $t_{12}\in \mathcal {T}_{12}$, using Lemma \ref{12}, we get
\begin{eqnarray*}
& &\phi (ct_{12})=\phi (c)\phi (t_{12})\\
&=&(\phi (a_{11}+\phi (b_{11}))\phi (t_{12})\\
&=&\phi (a_{11}t_{12})+\phi (b_{11}t_{12})\\
&=&\phi (a_{11}t_{12})+\phi (b_{11}t_{12}),
\end{eqnarray*}
which leads to $ct_{12}=a_{11}t_{12}+b_{11}t_{12}$. Accordingly,
$c_{11}=a_{11}+b_{11}$.
\end{proof}

\textbf{Proof of Theorem \ref{theorem}:} Suppose that $a$ and $b$
are two arbitrary elements of $\mathcal {T}$. We choose an element
$c\in \mathcal {T}$ satisfying $\phi (c)=\phi (a)+\phi (b)$. For any $t_{11}\in
\mathcal {T}_{11}$ and $s_{22}\in \mathcal {T}_{22}$, we obtain
\begin{eqnarray*}
\phi (t_{11}cs_{22})&=&\phi (t_{11})\phi (c)\phi (s_{22})\\
&=&\phi (t_{11})(\phi (a)+\phi (b))\phi (s_{22})\\
&=&\phi (t_{11})\phi (a)\phi (s_{22})+\phi (t_{11})\phi (b)\phi (s_{22})\\
&=&\phi (t_{11}as_{22})+\phi (t_{11}bs_{22})\\
&=&\phi (t_{11}as_{22}+t_{11}bs_{22}).
\end{eqnarray*}
Consequently, $c_{12}=a_{12}+b_{12}$.

Since $\phi $ is additive on $\mathcal {T}_{11}$, we can get
$c_{11}=a_{11}+b_{11}$ from $\phi (ct_{11})=\phi (at_{11})+\phi
(bt_{11})=\phi (at_{11}+bt_{11})$.

In the similar manner, one can get $c_{22}=a_{22}+b_{22}$. The proof
is complete.

 If algebras $\mathcal {A}$ and $\mathcal {B}$ contain identities, then we have the following result.
 \begin{corollary}
 Let $\mathcal {A}$ and $\mathcal {B}$ be two unital algebras over a commutative ring $\mathcal {R}$, $\mathcal {M}$ a faithful $(\mathcal {A}, \mathcal {B})$-bimodule, and $\mathcal {T}$ be the triangular algebra $Tri(\mathcal {A}, \mathcal {M}, \mathcal {B})$. Then any multiplicative bijective map from $\mathcal {R}$ onto an arbitrary ring $\mathcal {R}^{\prime }$ is additive.
 \end{corollary}

We end this section with the case when $\mathcal {A}$ and $\mathcal
{B}$ are standard operator algebras.
\begin{corollary}\label{coro1}
 Let $\mathcal {A}$ and $\mathcal {B}$ be two  standard operator algebras over a Banach space $X$, $\mathcal {M}$ a faithful $(\mathcal {A}, \mathcal {B})$-bimodule, and $\mathcal {T}$ be the triangular algebra $Tri(\mathcal {A}, \mathcal {M}, \mathcal {B})$. Then any multiplicative bijective map from $\mathcal {R}$ onto an arbitrary ring $\mathcal {R}^{\prime }$ is additive.
 \end{corollary}

\section{Jordan Maps on Triangular Algebras}

 In this section we deal with Jordan maps on triangular algebras.

Throughout this section, $\mathcal {T}=Tri(\mathcal {A}, \mathcal
{M}, \mathcal {B})$ will be a triangular algebra, where  $\mathcal
{A}$, $\mathcal {B}$ are two algebras over a commutative ring
$\mathcal {R}$ and $\mathcal {M}$ is a faithful $(\mathcal {A},
\mathcal {B})$-bimodule satisfying

(i) If $a\in \mathcal {A}$ and  $ax+xa=0$ for all $x\in \mathcal
{A}$, then $a=0$;

(ii)  If $b\in \mathcal {B}$ and  $by+yb=0$ for all $y\in \mathcal
{B}$, then $b=0$.

Map $\phi $ is a Jordan bijective map from $\mathcal {T}$ onto an
arbitrary ring $\mathcal {R}^{\prime }$.

We begin with the following lemma.
\begin{lemma}\label{0}
$\phi (0)=0$\end{lemma}
\begin{proof}
Pick $a\in \mathcal {T}$ such that $\phi (a)=0$. Then $\phi (0)=\phi
(a\cdot 0+0\cdot a)=\phi (a)\phi (0)+\phi (0)\phi (a)=0$.
\end{proof}

\begin{lemma} \label{2add}Suppose that $a, b, c\in \mathcal {T}$ satisfying $\phi (c)=\phi (a)+\phi (b)$, then for any $t\in \mathcal {T}$
$$\phi (tc+ct)=\phi (t)\phi (c)+\phi (c)\phi (t).$$
\end{lemma}
\begin{proof}
Multiplying $\phi (c)=\phi (a)+\phi (b)$ by $\phi (t)$ from the left
and the right respectively and adding them together, one can easily
get  $ \phi (tc+ct)=\phi (t)\phi (c)+\phi (c)\phi (t)$.
\end{proof}

\begin{lemma}\label{21112}
For any $a_{11}\in \mathcal {T}_{11}$ and $b_{12}\in \mathcal
{T}_{12}$, we have
$$\phi (a_{11}+b_{12})=\phi (a_{11})+\phi (b_{12}).$$
\end{lemma}
\begin{proof}
Let $c\in \mathcal {T}$ be chosen such that $\phi (c)=\phi
(a_{11})+\phi (b_{12})$. Now for any $t_{22}\in \mathcal {T}_{22}$,
by Lemma \ref{2add}, we have
$$\phi (t_{22}c+ct_{22})=\phi (t_{22}a_{11}+a_{11}t_{22})+\phi (t_{22}b_{12}+b_{12}t_{22})=\phi (b_{12}t_{22}).$$
It follows that $t_{22}c+ct_{22}=b_{12}t_{22}$, i.e.,
$t_{22}c_{22}+c_{12}t_{22}+c_{22}t_{22}=b_{12}t_{22}$. This implies
that $c_{12}t_{22}=b_{12}t_{22}$ and $t_{22}c_{22}+c_{22}t_{22}=0$,
and so $c_{12}=b_{12}$ and $c_{22}=0$.

From $$\phi (t_{12}c+ct_{12})=\phi (t_{12}a_{11}+a_{11}t_{12})+\phi
(t_{12}b_{12}+b_{12}t_{12})=\phi (a_{11}t_{12}),$$ one can get
$c_{11}=a_{11}$.
\end{proof}
Similarly, we have
\begin{lemma}\label{21222}
For arbitrary $a_{12}\in \mathcal {T}_{12}$ and $b_{22}\in \mathcal
{T}_{22}$, the following is true.
$$\phi (a_{12}+b_{22})=\phi (a_{12})+\phi (b_{22}).$$
\end{lemma}

\begin{lemma}\label{211121222}
$$\phi (a_{11}b_{12}+c_{12}d_{22})=\phi (a_{11}b_{12})+\phi (c_{12}d_{22})$$ holds true for any $a_{11}\in \mathcal {T}_{11}$, $b_{12}, c_{12}\in \mathcal {T}_{12}$, and $d_{22}\in \mathcal {T}_{22}$.
\end{lemma}
\begin{proof}
By Lemma \ref{21112} and Lemma \ref {21222}, we compute
\begin{eqnarray*}
& &\phi (a_{11}b_{12}+c_{12}d_{22})\\
&=&\phi ((a_{11}+c_{12})(b_{12}+d_{22})+(b_{12}+d_{22})(a_{11}+c_{12}))\\
&=&\phi (a_{11}+c_{12})\phi (b_{12}+d_{22})+\phi (b_{12}+d_{22})\phi (a_{11}+c_{12})\\
&=&(\phi (a_{11})+\phi (c_{12}))(\phi (b_{12})+\phi (d_{22}))+(\phi (b_{12})+\phi (d_{22}))(\phi (a_{11})+\phi (c_{12}))\\
&=&\phi (a_{11}b_{12}+b_{12}a_{11})+\phi (a_{11}d_{22}+d_{22}a_{11})+\phi (c_{12}b_{12}+b_{12}c_{12})+\phi (c_{12}d_{22}+d_{22}c_{12})\\
&=&\phi (a_{11}b_{12})+\phi (c_{12}d_{22}).
\end{eqnarray*}
\end{proof}
\begin{lemma}\label{212}
$\phi $ is additive on $\mathcal {T}_{12}$.
\end{lemma}
\begin{proof}
Let $a_{12}$ and $b_{12}$ be any two elements of $\mathcal
{T}_{12}$. Since $\phi $ is surjective, there exists a $c\in
\mathcal {T}$ such that $\phi (c)=\phi (a_{12})+\phi (b_{12})$.

Now for any $t_{22}\in \mathcal {T}_{22}$, by Lemma \ref{2add}, we
obtain
$$\phi (t_{22}c+ct_{22})=\phi (t_{22}a_{12}+a_{12}t_{22})+\phi (t_{22}b_{12}+b_{12}t_{22})=\phi (a_{12}t_{22})+\phi (b_{12}t_{22}).$$
Again, using Lemma \ref{2add}, for any $s_{11}\in \mathcal
{T}_{11}$, we have
\begin{eqnarray*}
& &\phi (s_{11}(t_{22}c+ct_{22})+(t_{22}c+ct_{22})s_{11})\\
&=&\phi (s_{11}a_{12}t_{22}+a_{12}t_{22}s_{11})+\phi (s_{11}b_{12}t_{22}+b_{12}t_{22}s_{11})\\
&=&\phi (s_{11}a_{12}t_{22})+\phi (s_{11}b_{12}t_{22})\\
&=&\phi (s_{11}a_{12}t_{22}+s_{11}b_{12}t_{22}).
\end{eqnarray*}
In the last equality we apply Lemma \ref{211121222}. It follows that
$$s_{11}c_{12}t_{22}=s_{11}a_{12}t_{22}+s_{11}b_{12}t_{22}.$$ This
gives us $c_{12}=a_{12}+b_{12}$.

To show $c_{11}=0$, we first consider $\phi (t_{11}c+ct_{11})$ for
any $t_{11}\in \mathcal {T}_{11}$. We have
$$\phi (t_{11}c+ct_{11})=\phi (t_{11}a_{12}+a_{12}t_{11})+\phi (t_{11}b_{12}+b_{12}t_{11})=\phi (t_{11}a_{12})+\phi (t_{11}b_{12}).$$

Furthermore, for arbitrary $s_{12}\in \mathcal {T}_{12}$,
$$\phi (s_{12}(t_{11}c+ct_{11})+(t_{11}c+ct_{11})s_{12})=\phi (s_{12}t_{11}a_{12}+t_{11}a_{12}s_{12})+\phi (s_{12}t_{11}b_{12}+t_{11}b_{12}s_{12})=0.$$
This implies that $t_{11}c_{11}s_{12}+c_{11}t_{11}s_{12}=0$, and so
$c_{11}=0$.

Note that $\phi (t_{12}c+ct_{12})=\phi
(t_{12}a_{12}+a_{12}t_{12})+\phi (t_{12}b_{12}+b_{12}t_{12})=0$.
Now,  $c_{22}=0$ follows easily.
\end{proof}
\begin{lemma}\label{211}
$\phi $ is additive on $\mathcal {T}_{11}$.
\end{lemma}
\begin{proof}
Suppose that $a_{11}$ and $b_{11}$ are two arbitrary elements of
$\mathcal {T}_{11}$. Let $c\in \mathcal {T}$ be an element of
$\mathcal {T}$ such that $\phi (c)=\phi (a_{11})+\phi (b_{11})$.

For any $t_{22}\in \mathcal {T}_{22}$, we get
$$\phi (t_{22}c+ct_{22})=\phi (t_{22}a_{11}+a_{11}t_{22})+\phi (t_{12}b_{11}+b_{11}t_{22})=0.$$
Therefore, $t_{22}c+ct_{22}=0$, which leads to $c_{12}=c_{22}=0$.

Similarly, we can get $c_{11}=a_{11}+b_{11}$ from
\begin{eqnarray*}
\phi (t_{12}c+ct_{12})&=&\phi (t_{12}a_{11}+a_{11}t_{12})+\phi (t_{12}b_{11}+b_{11}t_{12})\\
&=&\phi (a_{11}t_{12})+\phi (b_{11}t_{12})=\phi
(a_{11}t_{12}+b_{11}t_{12}).
\end{eqnarray*}
\end{proof}
\begin{lemma}\label{222}
$\phi $ is additive on $\mathcal {T}_{22}$.
\end{lemma}
\begin{proof}
For any $a_{22}, b_{22}\in \mathcal {T}_{22}$, by the surjectivity
of $\phi $, there is $c\in \mathcal {T}$ satisfying $\phi (c)=\phi
(a_{22})+\phi (b_{22})$.

We only show $c_{12}=0$. One can easily derive $c_{11}=0$ and
$c_{22}=a_{22}+b_{22}$ by considering $\phi (t_{11}c+ct_{11})$ and
$\phi (t_{12}c+ct_{12})$ for any $t_{11}\in \mathcal {T}_{11}$ and
$t_{12}\in \mathcal {T}_{12}$ respectively.

To show $c_{12}=0$. We first consider $\phi (t_{22}c+ct_{22})$ for
arbitrary $t_{22}\in \mathcal {T}_{22}$. We obtain
$$\phi (t_{22}c+ct_{22})=\phi (t_{22}a_{22}+a_{22}t_{22})+\phi (t_{22}b_{22}+b_{22}t_{22}).$$

For any $s_{11}$, using Lemma \ref{2add}, we get
\begin{eqnarray*}
& &\phi (s_{11}(t_{22}c+ct_{22})+(t_{22}c+ct_{22})s_{11})\\
&=&\phi (s_{11}(t_{22}a_{22}+a_{22}t_{22})+(t_{22}a_{22}+a_{22}t_{22})s_{11})\\& &+\phi (s_{11}(t_{22}b_{22}+b_{22}t_{22})+(t_{22}b_{22}+b_{22}t_{22})s_{11})\\
&=&0.
\end{eqnarray*}
Now, $c_{12}=0$ follows obviously.
\end{proof}

\begin{lemma}\label{21122}
For each $a_{11}\in \mathcal {T}_{11}$ and $b_{22}\in \mathcal
{T}_{22}$, we have
$$\phi (a_{11}+b_{22})=\phi (a_{11})+\phi (b_{22}).$$
\end{lemma}
\begin{proof}
Since $\phi $ is surjective, we can pick $c\in \mathcal {T}$ such
that $\phi (c)=\phi (a_{11})+\phi (b_{22})$.

Considering $\phi (t_{11}c+ct_{11})$ and $\phi (t_{22}c+ct_{22})$
for arbitrary $t_{11}\in \mathcal {T}_{11}$ and $t_{22}\in \mathcal
{T}_{22}$, one can infer that $c_{11}=a_{11}$, $c_{12}=0$, and
$c_{22}=b_{22}$.
\end{proof}

\begin{lemma}\label{2111222}
For any $a_{11}\in \mathcal {T}_{11}$, $b_{12}\in \mathcal
{T}_{12}$, and $c_{22}\in \mathcal {T}_{22}$,
$$\phi (a_{11}+b_{12}+c_{22})=\phi (a_{11})+\phi (b_{12})+\phi (c_{22}).$$
\end{lemma}
\begin{proof}
Let $d\in \mathcal {T}$ be chosen such that $\phi (d)=\phi
(a_{11})+\phi (b_{12})+\phi (c_{22})$. On one side, by Lemma
\ref{21112}, we have $$\phi (d)=\phi (a_{11}+b_{12})+\phi
(c_{22}).$$

Now for any $t_{11}\in \mathcal {T}_{11}$, we obtain
\begin{eqnarray*}
& &\phi (t_{11}d+dt_{11})\\
&=&\phi (t_{11}(a_{11}+b_{12})+(a_{11}+b_{12})t_{11})+\phi (t_{11}c_{22}+c_{22}t_{11})\\
&=&\phi (t_{11}a_{11}+t_{11}b_{12}+a_{11}t_{11}),
\end{eqnarray*}
 which gives us
$$t_{11}d_{11}+t_{11}d_{12}+d_{11}t_{11}=t_{11}a_{11}+t_{11}b_{12}+a_{11}t_{11}.$$
Hence, $d_{11}=a_{11}$ and $d_{12}=b_{12}$.

On the other side, by Lemma \ref{21122}, we see that
$$\phi (d)=\phi (a_{11}+c_{22})+\phi (b_{12}).$$
For any $t_{12}\in \mathcal {T}_{12}$, we have
\begin{eqnarray*}
& &\phi (t_{12}d+dt_{12})\\
&=&\phi (t_{12}(a_{11}+c_{22})+(a_{11}+c_{22})t_{12})+\phi (t_{12}c_{22}+c_{22}t_{12})\\
&=&\phi (t_{12}c_{22}+a_{11}t_{12}).
\end{eqnarray*}
We can infer $d_{22}=c_{22}$ from the fact that
$t_{12}d+dt_{12}=t_{12}c_{22}+a_{11}t_{12}$.
\end{proof}
We are in a position  to prove the main result of this section.

\begin{theorem}\label{2theorem}
Let $\mathcal {T}=Tri(\mathcal {A}, \mathcal {M}, \mathcal {B})$ be
a triangular algebra, where  $\mathcal {A}$, $\mathcal {B}$ are two
algebras over a commutative ring $\mathcal {R}$ and $\mathcal {M}$
is a faithful $(\mathcal {A}, \mathcal {B})$-bimodule. Suppose that
$\mathcal {A}$ and $\mathcal {B}$  satisfy

(i) If $a\in \mathcal {A}$ and  $ax+xa=0$ for all $x\in \mathcal
{A}$, then $a=0$;

(ii)  If $b\in \mathcal {B}$ and  $by+yb=0$ for all $y\in \mathcal
{B}$, then $b=0$.

Let $\phi $ be  a Jordan   map from $\mathcal {T}$ to an arbitrary
ring $\mathcal {R}^{\prime }$, i.e., for any $s, t\in \mathcal {T}$,
$$\phi (st+ts)=\phi (s)\phi (t)+\phi (t)\phi (s).$$
 If $\phi $ is bijective, then $\phi $ is additive.
\end{theorem}
\begin{proof}
For arbitrary $s$ and $t$ in $\mathcal {T}$. We write
$s=s_{11}+s_{12}+s_{22}$ and $t=t_{11}+t_{12}+t_{22}$. We compute
\begin{eqnarray*}
& &\phi (s+t)\\
&=&\phi ((s_{11}+s_{12}+s_{22})+(t_{11}+t_{12}+t_{22}))\\
&=&\phi ((s_{11}+t_{11})+(s_{12}+t_{12})+(s_{22}+t_{22}))\\
&=&\phi (s_{11}+t_{11})+\phi (s_{12}+t_{12})+\phi (s_{22}+t_{22})\\
&=&\phi (s_{11})+\phi (t_{11})+\phi (s_{12})+\phi (t_{12})+\phi (s_{22})+\phi (t_{22})\\
&=&(\phi (s_{11}+\phi (s_{12})+\phi (s_{22}))+(\phi t_{11}+\phi (t_{12})+\phi (t_{22}))\\
&=&\phi (s)+\phi (t).
\end{eqnarray*}The proof is complete.
\end{proof}

\section{Jordan Triple Maps on Triangular Algebras}

The aim of this section is to investigate the additivity of Jordan
triple maps on triangular algebras. We shall show that every Jordan
triple bijective map from a triangular algebra onto any ring is
automatically additive.  We only give the outline of the proof as it
is a modification of the proof of the related results of  section 2.
\begin{theorem}
  Let $\mathcal {A}$ and $\mathcal {B}$ be two algebras over a commutative ring $\mathcal {R}$, $\mathcal {M}$ be a faithful $(\mathcal {A}, \mathcal {B})$-bimodule, and $\mathcal {T}$ be the triangular algebra $Tri(\mathcal {A}, \mathcal {M}, \mathcal {B})$. Suppose that algebras $\mathcal {A}$ and $\mathcal {B}$ satisfy:

(i) If $a\mathcal {A}=\{ 0\} $, or $\mathcal {A}a=\{ 0\} $, then
$a=0$;

(ii) If $b\mathcal {B}=\{ 0\} $, or $\mathcal {B}b=\{ 0\} $, then
$b=0$.

Let $\psi $ be a Jordan triple map from $\mathcal {T}$ to an
arbitrary ring $\mathcal {R}^{\prime }$, i.e., for any $r, s, t\in
\mathcal {T}$,
$$\psi (rst+tsr)=\psi (r)\psi (s)\psi (t)+\psi (t)\psi (s)\psi (r).$$
If $\psi $ is bijective, then $\psi $ is additive.
\end{theorem}

\begin{proof} We divide the proof into a series of claims.
\begin{description}
\item[Claim 1] $\psi (0)=0$.

We find $a\in \mathcal{T}$ such that $\psi (a)=0$. Then $\psi
(0)=\psi (0\cdot a\cdot 0)+\psi (0\cdot a\cdot 0)=\psi (0)\psi
(a)\psi (0)+\psi (0)\psi (a)\psi (0)=0$.

\item[Claim 2]  If $\psi (c)=\psi (a)+\psi (b)$, then $$\psi (tsc+cst)=\psi (t)\psi (s)\psi (c)$$
holds true for all $s, t\in \mathcal {T}$.

One can get this easily by modifying the proof of Lemma \ref{2add}.

\item[Claim 3] $\psi (a_{11}+b_{12})=\psi (a_{11})+\psi (b_{12})$ and $\psi (a_{12}+b_{22})=\psi (a_{12})+\psi (b_{22})$.

We omit the proof as it is similar to the proof of Lemma \ref{21112}
and Lemma \ref{21222}.
\item [Claim 4]$\psi (t_{11}a_{11}b_{12}+t_{11}c_{12}d_{22})=\psi (t_{11}a_{11}b_{12})+\psi (t_{11}c_{12}d_{22})$.

We compute
\begin{eqnarray*}
& &\psi (t_{11}a_{11}b_{12}+t_{11}c_{12}d_{22})\\
&=&\psi (t_{11}(a_{11}+c_{12})(b_{12}+d_{22})+(b_{12}+d_{22})(a_{11}+c_{12})t_{11})\\
&=&\psi (t_{11})\psi (a_{11}+c_{12})\psi (b_{12}+d_{22})+\psi (b_{12}+d_{22})\psi (a_{11}+c_{12})\psi (t_{11})\\
&=&\psi (t_{11})(\psi (a_{11})+\psi (c_{12}))(\psi (b_{12})+\psi (d_{22}))\\
& &+(\psi (b_{12})+\psi (d_{22}))(\psi (a_{11})+\psi (c_{12}))\psi (t_{11})\\
&=&\psi (t_{11}a_{11}b_{12}+b_{12}a_{11}t_{11})+\psi (t_{11}a_{11}d_{22}+d_{22}a_{11}t_{11})\\
& &+\psi (t_{11}c_{12}b_{12}+b_{12}c_{12}t_{11})+\psi (t_{11}c_{12}d_{22}+d_{22}c_{12}t_{11})\\
&=&\psi (t_{11}a_{11}b_{12})+\psi (t_{11}c_{12}d_{22}).
\end{eqnarray*}
\item[Claim 5] $\psi $ is additive on $\mathcal {T}_{12}, \mathcal {T}_{11}$, and $\mathcal {T}_{22}$.

It follows similarly by modifying the proof of Lemma \ref{212},
Lemma \ref{211}, and Lemma \ref{222}, and using Claim 2, Claim 3,
and Claim 4.

\item [Claim 6] $\psi (a_{11}+b_{22})=\psi (a_{11})+\psi (b_{22})$.

See the proof Lemma \ref{21122}
\item [Claim 7] $\psi (a_{11}+b_{12}+c_{22})=\psi (a_{11})+\psi (b_{12})+\psi (c_{22})$.

See the proof of Lemma \ref{2111222}

\item[Claim 8] $\psi $ is additive.

The same as the proof of Theorem \ref{2theorem}

 \end{description}
\end{proof}

The following corollary follows directly if both $\mathcal {A}$ and
$\mathcal {B}$ are unital.
\begin{corollary}
Let $\mathcal {A}$ and $\mathcal {B}$ be two unital algebras over a
commutative ring $\mathcal {R}$, $\mathcal {M}$ be a faithful
$(\mathcal {A}, \mathcal {B})$-bimodule, and $\mathcal {T}$ be the
triangular algebra $Tri(\mathcal {A}, \mathcal {M}, \mathcal {B})$.
Suppose that $\psi $ is  a Jordan triple map from $\mathcal {T}$ to
an arbitrary ring $\mathcal {R}^{\prime }$. If $\psi $ is bijective,
then $\psi $ is additive.
\end{corollary}

Similar to Corollary \ref{coro1}, if both $\mathcal {A}$ and
$\mathcal {B}$ are standard operator algebras, we have
\begin{corollary}
Let $\mathcal {A}$ and $\mathcal {B}$ be two  standard operator
algebras over a Banach space $X$, $\mathcal {M}$ a faithful
$(\mathcal {A}, \mathcal {B})$-bimodule, and $\mathcal {T}$ be the
triangular algebra $Tri(\mathcal {A}, \mathcal {M}, \mathcal {B})$.
Then any Jordan triple bijective map from $\mathcal {R}$ onto an
arbitrary ring $\mathcal {R}^{\prime }$ is additive.
\end{corollary}

 \bibliographystyle{amsplain}

\begin{thebibliography}{10}



\bibitem{cheung} W. S. Cheung, Commuting maps on triangular algebras, \textit {J. London Math. Soc.}, \textbf {63} (2001), 117--127.

 \bibitem {ji} P. Ji,   Jordan maps on triangular algebras, \textit {Linear Algebra Appl.}, (to appear).

\bibitem {jing1} W. Jing, Aditivity of Jordan elementary maps on rings, preprint (2007).


\bibitem {jing2} W. Jing, Jordan triple elementary maps on rings, preprint (2007).

\bibitem {li237} P. Li, W. Jing, Jordan elementary maps on rings, \textit {Linear Algebra Appl.}, \textbf {382} (2004), 237--245.

\bibitem{lilu} P. Li, F. Lu, Additivity of elementary maps on rings, \textit{Comm. Algebra}, \textbf{32} (2004), 3725--3737.

\bibitem{lu123} F. Lu, Additivity of Jordan maps on standard operator algebras, \textit {Linear Algebra Appl.}, \textbf{357} (2002), 123--131.


\bibitem{lu2273} F. Lu, Jordan maps on associative algebras, \textit {Comm. Algebra}, \textbf {31} (2003), 2273--2286.
\bibitem{ma695} W. S. Martindale III, When are multiplicative mappings additive?\textit {Proc. Amer. Math. Soc.}, \textbf {21} (1969) 695--698.
 \end{thebibliography}

\end{document}